\documentclass[]{amsart}
\usepackage{graphicx}
\usepackage[T1]{fontenc}
\usepackage[latin1]{inputenc}
\usepackage{enumerate}
\usepackage{pstricks} 
\usepackage{latexsym,amssymb,amsthm,amsxtra}
\usepackage{amsmath}
\usepackage{stmaryrd}
\usepackage{tikz} 
\usepackage{pgf} 
{\bf}{\it}
\newtheorem{theorem}{Theorem}
\newtheorem{definition}{Definition}
\newtheorem{lemma}{Lemma}
\newtheorem{proposition}{Proposition}

\newtheoremstyle{hermesremark}{3pt}{3pt}{\normalfont}{0pt}{}{. --}{.5em}{}
\theoremstyle{hermesremark}
\newtheorem{remark}{{\bf Remark}}

\def\/{\, | \,}

\newcommand\prm{^{\prime}}

\def\llbracket{{\mathchoice {\rm |\mskip-4mu [} {\rm |\mskip-4mu [}
{\rm |\mskip-4mu [} {\rm |\mskip-4mu [}}}
\def\rrbracket{{\mathchoice {\rm ]\mskip-4mu |} {\rm ]\mskip-4mu |}
{\rm ]\mskip-4mu |} {\rm ]\mskip-4mu |}}}
\def\lescx{\prec_{\text{\textsc{scx}}}}
\def\lecx{\prec_{\text{\textsc{cx}}}}

\def\N{{\mathbb N}}

\newcommand{\R}{{\mathbb R}}

\def\ind{{\mathchoice {\rm 1\mskip-4mu l} {\rm 1\mskip-4mu l}
{\rm 1\mskip-4.5mu l} {\rm 1\mskip-5mu l}}}
\def\esp#1{{\mathbf E}\left[#1\right]}

\newcommand{\pae}[1]{\mbox{$\lfloor \kern-1pt #1 \kern-1pt \rfloor$}}
\newcommand{\paep}[1]{\mbox{$\lceil \kern-1pt #1 \kern-1pt \rceil$}}

\newcommand\suite[1]{\left\{#1_n;\,n\in\N\right\}}

\newcommand\suiten[1]{\left\{#1;\,n\in\N\right\}}

\def\pm{^\prime}

\newcommand\gre{\textbf{e}}

\begin{document}
\title{A pathwise comparison of parallel queues}
\author[P. Moyal]{Pascal Moyal}
\address{
Northwestern University\\
Department of Industrial Engineering \& Management Sciences\\
2145 Sheridan Road,\\
Evanston, IL 60208, USA}
\email{pascal.moyal@northwestern.edu}
\address{and}
\address{
Université de Technologie de Compiègne\\
Laboratoire de Mathématiques Appliquées de Compiègne\\
Rue Roger Couttolenc
CS 60319\\
60203 Compi\`egne, France.}
\email{pascal.moyal@utc.fr}

\begin{abstract}
 We introduce the appropriate framework for pathwise comparison of multiple server queues under general stationary ergodic assumptions. We show in what sense it is better to have more servers for a system under FCFS ('First Come, First Served') or equivalently, more queues in a system of parallel queues under the JSW ('Join the Shortest Workload') allocation policy. 
This comparison result is based on the recursive representation of Kiefer and Wolfowitz, and 
on a non-mass conservative generalization of the Schur-Convex semi-ordering.  
We also show that the latter result does not hold true in general, 
for the larger class of systems applying the semi-cyclic allocation policy introduced by Scheller-Wolf in \cite{SW03}.   
 
\end{abstract}
\maketitle

\section{Introduction}
\label{sec:intro}
Consider a queueing system of $S$  
{\em parallel queues}: there are $S$ servers, and to each one of them is associated a particular line. Moreoever the $S$ lines are independent of one another,   
so any incoming customer choses a server upon arrival and stays in the corresponding queue until service. 

The class of models we consider stipulate the knowledge of the {\em workload} ({\em i.e.} the quantity of work to be achieved) of each server, and apply the so-called {\em semi-cyclic service} (SCS) policy introduced by Scheller-Wolf in \cite{SW03}. 
Under an SCS policy, the server to which a customer is assigned 
is 'preserved', and won't welcome any other customer for the next $k$ arrivals, where $k \le S$ is fixed. Any incoming customer 
is then sent to the server having the least workload among the $S-k$ available ones. 
Particular cases of these policies are 'Join the Shortest Workload' (JSW), and the complete Cyclic case, respectively 
for $k=0$ and $k=S-1$.  
 
Another approach of queues with several servers assumes that the incoming customers are all put in the {\em same} queue during their wait. 
In the sequel, such models will be referred to as {\em multiple server queues}. % as opposed to the {\em parallel queues} which are studied in this paper. 
It is intuitively clear that multiple server queues are more flexible than parallel queues, in that they may offer to the customers the possibility of changing their service destination, to minimize their waiting time. In fact, parallel queues are nothing but a particular case of multiple server queues, having a service discipline which does not allow changes of queues. More precisely, it is easily seen that a parallel queue under the JSW policy is equivalent to a multiple 
server queue providing service in First Come, First Served (FCFS). 

For multiple server queues GI/GI, Foss \cite{Foss81,Foss89,CerFoss01} show in various distributional senses, the optimality of FCFS among a wide class of service disciplines (including those leading to parallel queues).  
This legitimates {\em in se} the use of the JSW policy, and this is why 
most of the literature on this topic focuses on this allocation policy. 

Under the most general assumption (stationary ergodic arrival process and 
service sequence), it leads to a simple representation by a stochastic recursive sequence, keeping track 
upon arrival times, of all residual service times, in increasing order 
(we call it the {\em service profile} of the system - see Kiefer and Wolfowitz equation (\ref{eq:KW})).  
 
Following this approach, Neveu \cite{Neu83} shows a stability condition for parallel queues under the JSW policy. This stability result is inherited from the monotonicity of the service profile sequence in some sense, and a minimal stationary profile is given by Loynes's Theorem \cite{Loynes62}. 
Foss \cite{Foss81} and Brandt \cite{Brandt85} then introduce the concept of maximal stationary profile for JSW queues. 

In this paper, we show that this representation provides the appropriate framework for a pathwise comparison of parallel queues. We show in what sense, a JSW system of $S$ queues performs better than one of $S\prm$ queues for $S\prm \le S$, providing 
both a transient result (Section \ref{sec:main}) and a comparison at equilibrium (Section \ref{sec:steady}). 
% Equivalently, we show {\em en route} that for a system of fixed size $S$, JSW performs better than a policy consisting in sending the arriving customers 
% to another (fixed) queue in the workload hierarchy. 
We concretely exhibit an ordering between vectors of different sizes ($S$ and $S\prm$) that is conserved recursively 
at any arrival time, whenever the two systems are coupled. The key technical points rely on the monotonicity of the 
recursive representation (\ref{eq:KW}), and on the introduction of 
a particular partial semi-ordering generalizing the Schur-convex ordering (\cite{Arnold87}, \cite{MarOl79}) to vectors of different total sums. 

Interestingly enough, we also show that the latter result does {\em not} hold true in general for SCS systems: we introduce  
a weaker ordering, that is conserved recursively between two SCS($k$) systems having respectively $S$ and $S\prm$ servers, 
but also show that this does not guarantee in general that the ordering between the {\em total workloads} of the two systems is
 conserved recursively - which is crucial for applications. 
We finally provide counter examples establishing that Cyclic systems of different sizes, and SCS systems with different values of $k$, are 
not comparable in this framework. 

Hence these results shade a different light on a good-sense observation: the more servers, the better. 
This fact is far from 
trivial when it comes to comparing the paths of the systems, and is true in this strong sense only in the very particular case 
of a JSW system of parallel queues, {\em i.e.} a FCFS multiple-server queue.

After we introduce precisely our model in section \ref{sec:model}, we present our main transient comparison results concerning SCS and JSW systems in section \ref{sec:main}, together with insightful counter-examples, and a discussion of the results in 
\cite{SW03}. Proofs of the latter results are provided in section \ref{sec:proof}. A natural extension of 
Theorem \ref{thm:mainJSW} to the comparison of JSW systems in steady state is proposed in section \ref{sec:steady}. 
We conclude this paper in section \ref{sec:comparealloc}, by giving insights on possible adaptations of the present approach 
to the comparison of allocation policies for systems of fixed size.

\section{The model}
\label{sec:model}

\subsection{Preliminary}
\label{sec:prelim}
\noindent
In what follows, $\R$ denotes the real line and $\R_+$, the subset of non-negative real numbers. Denote by $\N$ (respectively $\N^*$, $\mathbb Z$), the subset of non-negative (resp. positive, relative) integers. For any two elements $p$ and $q$ in $\N$, let 
$\llbracket p,q \rrbracket$ denote the collection  
$\{p,p+1,...,q\}$. For any $x,y \in \R$, denote $x \wedge y=\min\left(x,y\right)$ and $x \vee y=\max(x,y)$. Let $x^+=x \vee 0$. 

Fix $q \in \N^*$. We denote for all $x, y \in \left(\R_+\right)^q$ and $\lambda \in \R$, 
\begin{align*}
x&=\left(x(1),x(2),...,x(q)\right);\\ 
\lambda x&=\left(\lambda x(1),...,\lambda x(q)\right);\\
-x&=(-1)x\,;\\
x+y &=\left(x(1)+y(1),...,x(q)+y(q)\right);\\
\mathbf 0&=(0,\,...,\,0);\\
\mathbf 1&=(1,\,...,\,1);\\
\gre_i &=(0,\,...,\,\underbrace{1}_{i},\, ... ,0),\mbox{ for all }i\in \llbracket 1,q \rrbracket;\\
x^+ &=\Bigl(x(1)^+,x(2)^+,...,x(q)^+\Bigl);\\
\bar x,&\,\mbox{ the fully ordered version of }x,\mbox{ i.e. }\bar x(1) \le ... \le \bar x(p).
\end{align*}
For notational simplicity in several formulas, we adopt the convention that for any $x \in \left(\R_+\right)^q$, 
$x(q+1)=\infty$ . The coordinate-wise ordering on $\left(\R_+\right)^p$ is denoted "$\le$", {\em i.e.} $x \le y$ means that 
$x(1) \le y(1), ...,x(p) \le y(p)$.  
Let us also denote $\overline {(\R_+)^q}$, the subset of fully ordered vectors having non-negative coordinates.

\subsection{Systems of parallel queues}
\label{subsec:model}
Throughout, all random variables involved are defined on a common probability space $\mathcal Q:=\left(\Omega,\mathcal F, \mathbf P\right)$.  
We consider a queueing system having $S$ servers ($S \in \N^*$) working at the same pace. The customers enter the system according to a marked stationary point process, and we denote for all $n\in\N$, $\sigma_n$ the service time requested by customer $n$ and $\tau_n$, the $n$-th inter-arrival time, {\em i.e.} the duration between the arrivals of customers $n$ and $n+1$. Hence according to our assumptions the sequence $\suiten{\left(\sigma_n,\tau_n\right)}$ is identically distributed. 
We also assume that $\tau_0$ and $\sigma_0$ are both integrable, and that $\tau_0 >0$, $\mathbf P$-almost surely 
(a.s. for short), in other words the arrival point process is simple. 
Realizations of the sequence 
$\suiten{\left(\sigma_n,\tau_n\right)}$ will be called {\em input} of the system. 

We assume that the servers work in parallel, in the following sense: to each of the servers is associated a queue. 
The customers are allocated to a server upon arrival, and this choice (called hereafter {\em allocation policy}) is definitive. 
Inside each line, the corresponding server provides full service to its customers in First in, First out (FIFO). We assume that at each time, a full information is available on the {\em workload}, {\em i.e.} the quantity of work that each server still has to achieve, in time unit.

\subsubsection*{Semi-Cyclic systems} 
We now make precise the class of allocation policies we consider. Let $k \in \llbracket 0,S \rrbracket$. 
The {\em semi-cyclic service} policy SCS($k$) has been introduced in a general manner in \cite{SW03}. 
It consists in 'preserving' the servers to which the last $k$ entered customers have been affected, 
and sending any incoming customer to the server of shortest workload among the $S-k$ other ones. 
At any given time, the $k$ servers to which no new customer can be assigned are called {\em turn-over} servers (we also say that the turn-over is of size $k$).  In particular, if $k=S-1$ the policy is purely 'Cyclic' (Cy), in the sense that the incoming customers are assigned to the $S$ servers alternately, according to a strict round-robin rule. On another hand, if there is no turn-over ($k=0$), all servers are always available and the allocation policy is simply {\em Join the Shortest Workload} (JSW). As well known, the latter amounts to a multiple server queue where the $S$ servers operate in FCFS. 

Denote for any $n\ge 0$, $W_n^{k,S}=\left(W_n^{k,S}(1),...,W_n^{k,S}(S)\right)$ the {\em workload profile} seen by customer 
$n$ upon arrival, in a system of $S$ servers and a turn-over of $k$ servers. More precisely, for any $i \in \llbracket 1,S-(k\wedge n) \rrbracket$, $W_n^{k,S}(i)$ is the $i$-th workload in the increasing order among the 
$S-(k\wedge n)$ available servers, and for any $\ell\in \llbracket 0,(k\wedge n)-1 \rrbracket$,  
$W_n^{k,S}(S-\ell)$ denotes the workload of the server to which customer $n-(\ell+1)$ was assigned.  
The total {\em workload} $\mathbf W_n^{k,S}$ seen by customer 
$n$ upon arrival is then given by the sum of the workload profile, {\em i.e.} 
$$\mathbf W_n^{k,S} = \sum_{i=1}^S W_n^{k,S}(i),\quad n\ge 1.$$  
Then the service time $\sigma_n$ brought by customer $n$ is assigned to the server of workload $W_n^{k,S}(1)$. 

\subsubsection*{Recursive representation} Clearly, the sequence $\suite{W^{k,S}}$ is stochastically recursive, in the sense 
that for all $n$, its value at rank $n+1$ is fully determined by that at rank $n$ and the $n$-th component of the input. 
Specifically, 

\medskip

\noindent $\underline{\mbox{Case }k>0}$. 
When there is a turn-over, we have almost surely for all $n\ge 0$,   
\begin{equation*}
%\label{eq:recurSCSkS1}
\left\{\begin{array}{ll}
\Bigl(W_{n+1}^{k,S}(1),...,W_{n+1}^{k,S}(S-k)\Bigl)&=\overline{\biggl(\Bigl[W_{n}^{k,S}(2)-\tau_n\Bigl]^+,...,\Bigl[W_{n}^{k,S}(S-k+1)-\tau_n\Bigl]^+\biggl)};\\
							&\\
W_{n+1}^{k,S}(i) &= \left[W_{n}^{k,S}(i+1)-\tau_n\right]^+, i\in \llbracket S-k+1,S-1 \rrbracket;\\
							&\\
W_{n+1}^{k,S}(S) &= \left[W_{n}^{k,S}(1)+\sigma_n-\tau_n\right]^+.
							\end{array}\right.							
\end{equation*}
Equivalently, we have that 
\begin{equation}\label{eq:recurSCSkS1}
W_{n+1}^{k,S}=\Phi^{k,S}\left(W_{n}^{k,S},\sigma_n,\tau_n\right),\,n\ge 0,\,\mbox{ a.s.,}
\end{equation}
where for all $s,t\ge 0$, all $u\in \left(\R_+\right)^S$ and any $i\in \llbracket 1,S \rrbracket$, 
\begin{equation}
\label{eq:recurSCSkS2}
\Phi^{k,S}(u,s,t)(i)=\left\{\begin{array}{ll}
                            &\left[u(2) \wedge u(S-k+1)-t\right]^+\mbox{ if }i=1;\\ 
							&\\
							&\left[u(i) \vee \left(u(i+1) \wedge u(S-k+1)\right)-t\right]^+\mbox{ if }i\in \llbracket 2,S-k \rrbracket;\\ 
							&\\
							&\left[u(i+1)-t\right]^+\mbox{ if } i\in \llbracket S-k+1,S-1 \rrbracket;\\
							&\\
							&\left[u(1)+s-t\right]^+\mbox{ if }i=S.
							\end{array}\right.							
\end{equation}
We say in such a case that $\suite{W^{k,S}}$ is the $(\R_+)^S$-valued {\em Stochastic Recursive Sequence} (SRS) driven by the 
{\em driving map} $\Phi^{k,S}(.,.,.)$.

\medskip

\noindent $\underline{\mbox{Case }k=0}$. 
For the Join the Shortest Workload (or equivalently, FIFO) system, the driving map is given by the well-known Kiefer and Wolfowitz relation 
\cite{KW55}: a.s. for all $n \ge 0$,    
\begin{equation}
\label{eq:KW}
W_{n+1}^{0,S}=\overline{\biggl(\Bigl[W_{n}^{0,S}(1)+\sigma_n-\tau_n\Bigl]^+,\Bigl[W_{n}^{0,S}(2)-\tau_n\Bigl]^+,...,\Bigl[W_{n}^{0,S}(S)-\tau_n\Bigl]^+\biggl)}.
\end{equation}
%\label{eq:recurSCSkS1}
This amounts to writing that 
\begin{equation}\label{eq:recurSCS0S1}
W_{n+1}^{0,S}=\Phi^{0,S}\left(W_{n}^{0,S},\sigma_n,\tau_n\right),\,n\ge 0,\,\mbox{ a.s.,}
\end{equation}
where the driving map $\Phi^{0,S}(.,.,.)$ is defined, for all $s,t\ge 0$ and all $u\in \left(\R_+\right)^S$ by
\begin{equation}
\label{eq:recurSCS0S2}
\Phi^{0,S}(u,s,t)(i)=\Bigl[u(i) \vee \left(\left(u(1)+s\right)\wedge u(i+1)\right)-t\Bigl]^+,\,i\in \llbracket 1,S \rrbracket.
\end{equation}

\subsubsection*{Join the $p$-th shortest workload} 
Let $p \le S$. In the proofs below, we will also be led to consider a system of $S$ servers, without turn-over, in which the incoming customers are always sent to the server having the $p$-th smallest workload. The workload profile sequence $\suite{W^{0,S,p}}$ of this model is clearly stochastic recursive: we have 
\begin{equation}
\label{eq:recurJpSW}
W_{n+1}^{0,S,p}=\Phi^{0,S,p}\left(W_{n}^{0,S,p},\sigma_n,\tau_n\right),\,n\ge 0,\,\mbox{ a.s.}
\end{equation}
where, similarly to (\ref{eq:recurSCS0S2}), the driving map $\Phi^{0,S,p}(.,.,.)$ is given 
for all $u \in \left(\R_+\right)^S$ and $s,t \ge 0$ by  
\begin{equation}
\label{eq:recurSCS0SP2}
\Phi^{0,S,p}(u,s,t)(i)=\Bigl[u(i) \vee \left(\Bigl(u(p)+s\Bigl)\wedge u(i+1)\right)-t\Bigl]^+,\,i\in \llbracket 1,S \rrbracket.
\end{equation} 
Observe that the latter allocation policy does not make much sense in practice, since it leaves at least $p-1$ servers eventually inactive (depending on the 
initial conditions). We elaborate in Section \ref{sec:comparealloc}. In fact, as will appear clearly in the proof of Theorem \ref{thm:mainJSW}, 
this alternate policy will essentially be used as a tool for comparing JSW systems of different sizes.

\section{Main results}
\label{sec:main}
We establish in this paper in what sense, increasing the number of servers decreases, 
{\em on a pathwise basis}, the workload of semi-cyclic systems. Let us first introduce the following binary relations 
"$\prec$" and "$\prec_*$" on $\left(\R_+\right)^S \times \left(\R_+\right)^{S\pm}$, where $S\pm$ and $S$ are two positive 
integers such that $S\pm \le S$ (we omit the dependence in $S,S\pm$ for notational simplicity). 
\begin{definition}
\label{def:ordreétoile}
Let $S\prm \le S$ and let $u\in \left(\R_+\right)^S$ and $v \in \left(\R_+\right)^{S\pm}$.  
\begin{itemize}
% \item[(i)] We denote $u \prec_{\text{\tiny{wcx}}} v$ if 
% $$\sum_{i = k}^S \overline u(i) \le \sum_{i = k}^{S\pm} \overline v(i)\quad\mbox{ for any }k\in \llbracket 1,S \rrbracket.$$
\item[(i)] We write $u \prec v$ if
\[\left\{\begin{array}{ll}
         & u(S-\ell) \le v(S\pm-\ell)\,\mbox{ for all }\ell \in \llbracket 0,S\pm-1 \rrbracket;\\
         & \displaystyle\sum_{i=1}^S u(i) \le \displaystyle\sum_{i=1}^{S\pm} v(i).
         \end{array}\right.\]
\item[(ii)] We write $u \prec_{*} v$ if
\[\left\{\begin{array}{ll}
         & u(S-\ell) \le v(S\pm-\ell)\,\mbox{ for all }\ell \in \llbracket 0,S\pm-2 \rrbracket;\\
         & u(1) \le v(1).
         \end{array}\right.\]
\end{itemize} 
\end{definition}

\medskip

The binary relation "$\prec$" on $\left(\R_+\right)^S \times \left(\R_+\right)^{S\pm}$ can be interpreted as the most natural 
generalization of the coordinate-wise ordering "$\le$" to vector of different sizes: $u$ is less than $v$ in this sense 
if the restriction of $u$ to its last $S\pm$ coordinates is less than $v$ for "$\le$" on $\left(\R_+\right)^{S\pm}$, and if furthermore the {\em sums} of all coordinates of $v$ and $u$ and ordered in the same sense. 
% The relation "$\prec_*$" is a slight variant of $\prec$, weaker 
% in the present context, for $u \prec v$ entails $u \prec_* v$ whenever $u(S-S\pm+1) \ge u(1)$. 

The latter orderings "$\prec$" and "$\prec_*$"  turn out to be the suitable tool for comparing systems of parallel queues of different sizes. We have the following results, 

\begin{theorem}[Semi-cyclic systems]
\label{thm:mainSCS}
Let $2\le S \le S\pm$, and let $k$ be an integer such that $0 \le k \le S\pm-2$. 
Consider two systems of type SCS ($k$), having respectively $S$ and $S\prm$ servers, and fed 
by the same input $\suiten{\left(\sigma_n,\tau_n\right)}$. 
Assume that the corresponding original service profiles satisfy $W_0^{k,S} \prec_* W_0^{k,S^\prime}$. 
Then, this is true at all times: for all $n\ge 0$, $W_n^{k,S} \prec_* W_n^{k,S^\prime}$, or in other words 
\begin{align}
W_n^{k,S}(S-i) &\le W_n^{k,S\prm}(S\prm-i),\,i\in \llbracket 0,S\prm-2\rrbracket;\label{eq:main1SCS}\\
W_n^{k,S}(1) &\le W_n^{k,S\prm}(1).\label{eq:main1SCS}
% \mathbf W_n^{k,S} &\le \mathbf W_n^{k\prm,S\prm}.\label{eq:main2SCS}
\end{align}
\end{theorem}

\medskip

\begin{theorem}[JSW systems]
\label{thm:mainJSW}
Consider two SCS($0$) ({\em i.e.} JSW) systems, having respectively $S$ and $S\pm$ servers, 
and fed by the same input $\suiten{\left(\sigma_n,\tau_n\right)}$. 
Assume that the corresponding original service profiles are such that $W_0^{0,S} \prec W_0^{0,S\pm}$. 
% \begin{align*}
% W_0^{0,S}(S-i) &\le W_0^{0,S\prm}(S\pm-i),\,i\in \llbracket 0,S\prm-1\rrbracket;\\
% \mathbf W_0^{0,S} &\le \mathbf W_0^{0,S\prm}.
% \end{align*} 
Then, this is true at all times: for all $n\ge 0$, $W_n^{0,S} \prec W_n^{0,S\pm}$ or in other words, 
\begin{align}
W_n^{0,S}(S-i) &\le W_n^{0,S\prm}(S\pm-i),\,i\in \llbracket 0,S\prm-1\rrbracket;\label{eq:main1}\\
\mathbf W_n^{0,S} &\le \mathbf W_n^{0,S\prm}.\label{eq:main2}
\end{align}
\end{theorem}

\medskip

% \bigskip
% Finally, we show that Theorem \ref{thm:mainJSW} can be extended to the stady state of the systems under consideration, 
% as follows. 
% \begin{corollary}
% \label{cor:steady}
% Under the assumptions of Theorem \ref{thm:main}, assume that $\esp{\sigma} < S\prm\esp{\tau}.$ 
% Also assume that the sequence $\suiten{\left(\sigma_n,\tau_n\right)}$ is i.i.d.. 
% Let $W^{k,S}$ (respectively, $W^{k\prm,S\prm}$) be a $\left(\R_+\right)^S$-valued (resp., $\left(\R_+\right)^{S\prm}$-valued) random variable distributed 
% according to the unique stationary distribution of the Markov chain $\suiten{W^{k,S}_n}$ (resp., $\suiten{W^{k,S}}_n$). Then, we have 
% \begin{align*}
% W^{k,S}(S-i) &\le_{st} W^{k\prm,S\prm}(S\prm-i),\,i\in \llbracket 0,S\prm-1\rrbracket;\\
% \mathbf W^{k,S} &\le_{st} \mathbf W^{k\prm,S\prm},
% \end{align*}
% where "$\le_{st}$" denotes the strong stochastic ordering. 
% \end{corollary}

\medskip
Theorems \ref{thm:mainSCS} and \ref{thm:mainJSW} thus make precise in what sense an increase in the number of servers increases the efficiency 
of the system, on a pathwise basis and not only for some distributional order, respectively for SCS and JSW systems. 
 
A crucial observation resides in what Theorem \ref{thm:mainSCS} actually {\em does not} claim. First observe that the order 
"$\prec_*$" that is preserved by the recursion $\suite{W^{k,.}}$ is weaker than "$\prec$". In fact, $u \prec v$ implies $u \prec_* v$ for any $v$ whenever $u$ is ordered between coordinates 1 and $S-S\prm+1$ (which is the case for the workload profiles of a SCS system), but 
the converse is false. In particular, $u \prec_* v$ does not entail that $\sum_{i=1}^{S} u(i) \le \sum_{i=1}^{S\prm} v(i)$. In other words in this context, for a SCS($k$) system with $k>0$, even though the {\em total} workload of the small system is larger than that of the big system, the total workloads after one iteration need not be ordered in the same way. %We give an example below. 
We give an example below. 

For that reason, JSW systems ({\em i.e.} FIFO many-server systems) play a very singular role in this pathwise analysis, 
as JSW is the only SCS policy preserving the ordering of total workloads. This fact has a clear meaning in practical cases: minimizing the workload of the individual servers by increasing the number of these servers, is vain in terms of quality of service, if it does not lead to also minimize the total workload. 
%As is proven in this paper, it does so only for JSW systems, .  
% We draw hereafter a list of counter-examples to emphasize that such strong a result as Theorem \ref{thm:mainJSW} does not 
% hold for any other semi-cyclic system. 

\medskip

\subsection{Counter examples} 
\label{subsec:counter}
Let us provide several counter examples emphasizing that Theorems  
\ref{thm:mainSCS} and \ref{thm:mainJSW} cannot be generalized for a larger class of semi-cyclic systems: in \ref{subsubsec:cyclic} and 
\ref{subsubsec:diffk}, we show that Theorem \ref{thm:mainSCS} does not hold true for Cyclic systems and SCS systems with different 
turn-over sizes. In \ref{subsubsec:samek}, we show that the assumptions of Theorem \ref{thm:mainJSW} cannot be relaxed to include 
SCS systems with the same turn-over size. 

\subsubsection{Purely cyclic systems.}\label{subsubsec:cyclic}
 Let us first consider Cyclic systems, {\em i.e.} systems of $S$ servers with 
a turn-over of size $S-1$. Take for instance $S=3,k=2$, $S\prm=2$ and $k\prm=1$.  
Let $W_0^{2,3}=(9,3,1)$ and $W_0^{1,2}=(3,11)$, in a way that $W_0^{2,3} \prec W_0^{1,2}$. 
Set $\left(\sigma_0,\tau_0\right)=(4,0.5)$. We then have $W_{1}^{2,3}=(2.5,0.5,12.5)$ and $W_{1}^{1,2}=(10.5,6.5)$, thus 
$W_1^{2,3} \not\prec W_1^{1,2}$. Then, if $\left(\sigma_{1},\tau_{1}\right)=(1,12)$, we have $W_{2}^{2,3}=(0,0,0.5)$ and $W_{2}^{1,2}=(0,0)$, in particular the total workload at time $2$ is bigger in the large system. 

\subsubsection{Semi-cyclic systems with different turn-over sizes.}\label{subsubsec:diffk}
Set $k=2,S=4$, $k\prm=1$ and $S\prm=3$. 
Set for instance $W_0^{2,4}=(3,6,7,2)$ and $W_0^{1,3}=(7,8,3)$, and take $(\sigma_0,\tau_0)=(3,1)$. 
The workload profiles are then updated as follows: 
$$W_{1}^{2,4}=(5,6,1,5);\quad W_{1}^{1,3}=(2,7,9).$$
We thus have $W_0^{2,4} \prec W_0^{1,3}$ and $W_0^{2,4} \prec_* W_0^{1,3}$,  
but $W_{1}^{2,4} \not\prec W_{1}^{1,3}$ and $W_{1}^{2,4} \not\prec_* W_{1}^{1,3}$.  
Further, if for instance 
$\left(\sigma_1,\tau_1\right)=(9,12)$ we have $W_{2}^{2,4}=(0,0,0,2)$ and $W_{2}^{1,3}=(0,0,0)$, so the total workload 
after 2 iterations is again larger in the large system.  

\subsubsection{Semi-cyclic systems with equal turn-over sizes.} 
\label{subsubsec:samek}
Theorem \ref{thm:mainJSW} does not hold true even for SCS systems having the same turn-over sizes. As a counterexample, take $S=4,k=1$, $S\prm=3$ and $k\prm=1$ and set  
$W_{0}^{1,3}=(5,6,3)$ and $W_{0}^{1,4}=(1,4,6,3)$, in a way that $W_0^{1,4} \prec W_0^{1,3}$ and 
$W_0^{1,4} \prec_* W_0^{1,3}$. Then for $(\sigma_0,\tau_0)=(4,1)$ we obtain that $W_{1}^{1,3}=(2,5,8)$ and $W_{1}^{1,4}=(2,3,5,4)$, so 
$W_{1}^{1,4} \prec_* W_{1}^{1,3}$ but $W_{1}^{1,4} \not\prec W_{1}^{1,3}$. 

Let us now emphasize that the ordering of sums is not conserved for SCS systems having the same turn-over size. 
Take $S=5,k=2$, $S\prm=4$ and $k\prm=2$ and set for instance $W_{0}^{2,4}=(6,12,10,10)$ and $W_{0}^{2,5}=(4,11,12,10,1)$, 
in a way that $W_0^{2,5} \prec_* W_0^{2,4}$ and $\mathbf W_0^{2,5} \le \mathbf W_0^{2,4}$. Then, if for example $\left(\sigma_0,\tau_0\right)=(4,10)$ we obtain $W_{1}^{2,4}=(0,2,0,0)$ and $W_{1}^{2,5}=(0,1,2,0,0)$. 
Therefore, even though $W_1^{2,5} \prec_* W_1^{2,4}$ we have $\mathbf W_1^{2,5} > \mathbf W_1^{2,4}$.

% \medskip

% It is interesting to observe that in all cases, the same holds true for 
% the ordered versions of the profiles. In fact we have 
% $W_0^{k,S} \prec_{S\prm}^{S} W_0^{k\prm,S\prm}\mbox{ and }\overline{W_0^{k,S}} \prec_{S\prm}^{S} \overline{W_0^{k\prm,S\prm}}$, 
% but  
% $W_{1}^{k,S} \not{\prec_{S\prm}^{S}} W_{1}^{k\prm,S\prm}\mbox{ and }\overline{W_{1}^{k,S}} \not{\prec_{S\prm}^{S}} \overline{W_0^{k\prm,S\prm}}.$ 
% Further, observe that $\mathbf W_{2}^{k,S} >  \mathbf W_{2}^{k\prm,S\prm}$ in all scenarios, in other words we end up with 
% a smaller total workload in the smaller system.  

\subsection{About Lemma 4.1 in \cite{SW03}}
\label{subsec:SW}
These results have strong connections with Lemma 4.1 of \cite{SW03} which states, in our settings, that for all $S$ and all $0\le k \le S-2$, 
whenever $W_0^{k,S}=\mathbf 0$ and $W_0^{k,S-1}=\mathbf 0$ and the two systems are fed by the same input, we have  
\begin{equation*}
%\label{eq:SW1}
W_n^{k,S}(S-i) \le W_n^{k,S-1}(S-1-i),\,i\in \llbracket 0,S-2\rrbracket,\,n\ge 0.
\end{equation*}
In fact, the argument of the proof of this lemma stems itself from an inductive argument, from which it readily follows that 
\begin{multline}
\label{eq:SW1}
\left[W_0^{k,S}(S-i) \le W_0^{k,S-1}(S-1-i),\,i\in \llbracket 0,S-2\rrbracket\right] \\
\Longrightarrow \left[W_n^{k,S}(S-i) \le W_n^{k,S-1}(S-1-i),\,i\in \llbracket 0,S-2\rrbracket,\,n\ge 0\right]. 
\end{multline}
In turn, (\ref{eq:SW1}) readily entails by induction on $S\prm$ that for all $S\prm \le S$ such that $k \le S-2$, 
\begin{multline}
\label{eq:SW2}
\left[W_0^{k,S}(S-i) \le W_0^{k,S\prm}(S\prm-i),\,i\in \llbracket 0,S\prm-1 \rrbracket\right] \\
\Longrightarrow \left[W_n^{k,S}(S-i) \le W_n^{k,S\prm}(S\prm-i),\,i\in \llbracket 0,S\prm-1 \rrbracket,\,n\ge 0\right]. 
\end{multline} 
% Therefore, Theorem \ref{thm:mainJSW} and \ref{thm:mainSCS} have to be related to (\ref{eq:SW2}) and thereby, 
% to Lemma 4.1 in [\emph{ibid.}]. 
We make the following remarks: 

\medskip

(1) The inductive argument in the proof of Lemma 4.1 in [\emph{ibid.}] (more precisely the second display in that proof), 
leading to (\ref{eq:SW1}) above, is false: we give a counter-example in sub-section \ref{subsubsec:samek}. 
What is true is only the weaker statement of Theorem \ref{thm:mainSCS}; 

\medskip

(2) As a matter of fact, Lemma 4.1 in [\emph{ibid.}] holds true only for $k=0$, {\em i.e.} for a JSW 
system (Theorem \ref{thm:mainJSW}). But the proof in [\emph{ibid.}] (and more precisely the inequality in the first display) does not hold for $k=0$ 
and supposes the existence of a turn-over. We provide a correct proof of the recursive property of the statement 
(\ref{eq:main1}) in the proof of Theorem \ref{thm:mainJSW}. 
% The crucial case of a FIFO many-server queue ({\em i.e.} $k=0$) is therefore not covered by Lemma 4.1 in [\emph{ibid.}], but is addressed in Theorem \ref{thm:main} above;

\medskip

(3) Last, we complete the comparison of JSW systems by the second inequality (\ref{eq:main2}), stating 
that the larger systems minimizes the total workload on a pathwise basis. This fact is of course, not implied by (\ref{eq:main1}). 

% \medskip

% The rest of the paper is devoted to the proofs of Theorem \ref{thm:mainSCS}, \ref{thm:mainJSW} and Corollary \ref{cor:steady}. 
% For the latter, a stronger result is stated in Section \ref{sec:steady} under general stationary ergodic assumptions.  

\section{Proofs}
\label{sec:proof}

\subsection{Proof of Theorem \ref{thm:mainSCS}}
\label{subsec:proofSCS}
We start by proving that the ordering "$\prec_*$" is preserved by the driving maps $\Phi^{k,S}$ and $\Phi^{k,S\pm}$ for any 
fixed $k$. 
\begin{proposition}
\label{pro:mainSCS}
Let $S,S\pm \ge 3$ such that $S \le S\prm$, and fix $k\in \llbracket 1,S\pm-2 \rrbracket$. Let $u \in (\R_+)^S$ and $v \in \left(\R_+\right)^{S\pm}$ such that 
$u \prec_* v$. Then, for all $s,t \in \R_+$ we have that 
$$\Phi^{k,S}(u,s,t) \prec_* \Phi^{k,S\pm}(v,s,t).$$ 
\end{proposition}
\begin{proof}
Suppose that $u \prec_* v$ and fix $s$ and $t$. 
First, observe that 
$$\Phi^{k,S}(u,s,t)(S)=\left[u(1)+s-t\right]^+ \le \left[v(1)+s-t\right]^+=\Phi^{k,S\pm}(v,s,t)(S\pm).$$
Then, for all $\ell \in \llbracket 1,k-1 \rrbracket$, 
$$\Phi^{k,S}(u,s,t)(S-\ell)=\left[u(S-\ell+1)-t\right]^+ \le \left[v(S\pm-\ell+1)-t\right]^+=\Phi^{k,S\pm}(v,s,t)(S\pm-\ell).$$
Now, take $\ell \in \llbracket k,S\pm-2 \rrbracket.$ We have that 
\begin{align*}
\Phi^{k,S}(u,s,t)(S-\ell) &= \left[u(S-\ell) \vee \left(u(S-\ell+1) \wedge u(S-k+1)\right)-t\right]^+\\
                                &\le \left[v(S\pm-\ell) \vee \left(v(S\pm-\ell+1) \wedge v(S\pm-k+1)\right)-t\right]^+\\
                                &= \Phi^{k,S\pm}(v,s,t)(S\pm-\ell).
                                \end{align*}
Only the index $\ell=S\pm-1$ remains to be treated. We have 
\begin{align*}
\Phi^{k,S}(u,s,t)(1)    &= \left[u(2) \wedge u(S-k+1)-t\right]^+\\
                                &\le \left[u(S-S\pm+2) \wedge u(S-k+1)-t\right]^+\\
                                &\le \left[v(2) \wedge v(S\pm-k+1)-t\right]^+\\
                                &= \Phi^{k,S\pm}(v,s,t)(1),
                                \end{align*}
which concludes the proof. 
\end{proof}

\begin{proof}[Proof of Theorem \ref{thm:mainSCS}]
The result follows from Proposition \ref{pro:mainSCS} and the recursive representation (\ref{eq:recurSCSkS1}), 
using a simple induction.  
\end{proof}

\subsection{Proof of Theorem \ref{thm:mainJSW}}
\label{subsec:proofJSW}
As we will show, the proof of Theorem \ref{thm:mainJSW}, 
which relies in particular on the fact that the ordering of sums is conserved by the 
maps $\Phi^{0,S}$ and $\Phi^{0,S\pm}$, has a flavor of convex comparison. Using this approach we in fact show a more general result, Proposition \ref{pro:mainJSW} below. Let us first recall the definition and two main properties of the Schur-convex ordering "$\lescx$". 

\begin{definition}[Schur-convex ordering]
\label{def:schurconvex}
Let $S\ge 1$. For all $u,v \in \left(\R_+\right)^S$, we write $u \lescx v$ whenever 
\[\left\{\begin{array}{ll}
\displaystyle\sum_{i=1}^S u(i)=\displaystyle\sum_{i=1}^S v(i),&\\
\displaystyle\sum_{i=k}^S \bar u(i) \le \displaystyle\sum_{i=k}^S \bar v(i),&\mbox{ for all }k \in \llbracket 2,S \rrbracket.
\end{array}\right.\]
\end{definition}
The reader is referred to \cite{Arnold87}, \cite{MarOl79} and the Chapter 4 of \cite{BacBre02}  for an exhaustive presentation of the 
properties of the partial semi-ordering "$\lescx$". Let us quote two particular ones, which will be used hereafter.   
First, for any $u \in \left(\R_+\right)^S$ and any permutation $\alpha$ of $\llbracket 1,S \rrbracket$, define 
$$u_\alpha=\left(u\left(\alpha(1)\right),...,v\left(\alpha(S)\right)\right).$$ 
Observe that reordering twice the coordinates of a vector $u$, first according to the permutation $\alpha$ and then 
according to $\beta$, yields to the vector $u_{\alpha\circ\beta}$. 

A mapping $F:\,\left(\R_+\right)^S \to E$ is termed \emph{symmetric} whenever $F(u)=F\left(u_{\alpha}\right)$  for all $u \in \left(\R_+\right)^S$ and all permutations $\alpha$. Finally, $\alpha$ is said to be a {\em reordering} permutation  
of the vector $u$ if for some couple of indexes $i < j$ such that $u(i) \ge u(j)$, we have $\alpha(i)=j$, $\alpha(j)=i$ 
and $\alpha(\ell)=\ell$ for any $\ell \not\in \{i,j\}$. 
The following result can be found \emph{e.g.} in \cite{BacBre02} (Prop. 4.1.1 p.262 and Lemma 4.1.1 p.266): 
\begin{lemma}
\label{lemma:schur-convexe1} 
\begin{enumerate}
\item[(i)] For all $u,v \in \left(\R_+\right)^S$,   
\begin{equation*} 
 u \lescx v \Longleftrightarrow 
\Bigl[F(u)\le F(v) \mbox{ for all convex symmetric mappings }F: \left(\R_+\right)^S \to \R\Bigl].
\end{equation*}
\item[(ii)]  For all $u,v \in \left(\R_+\right)^S$ and any reordering permutation $\alpha$ of $u$, we have 
$$u_{\alpha} - \bar v \lescx u - \bar v.$$ 
\end{enumerate}
\end{lemma}

\noindent We deduce the following result, which states in the queueing context, that sending a new job to a 
more loaded server is always less optimal, in the convex sense, than to a less loaded server.    
\begin{lemma}
\label{lemma:schur-convexe3}
Let $w$ be an element of $(\R_+)^{q}$, where $q\ge 2$. Let $i,j \in \llbracket 1,q \rrbracket$ such that 
$w(i) \le w(j)$. Then, for any two nonnegative real numbers $s,t$, 
$$w+s\gre_i-t\mathbf 1 \lescx w+s\gre_j-t\mathbf 1.$$
\end{lemma}

\begin{proof}
Let us define the following permutations $\alpha$ and $\beta$ of $\llbracket 1,q \rrbracket$: 
\begin{enumerate}
\item if $1 \not\in \{i,j\}$, then $\alpha$ is the cycle $(j,i,1)$, {\em i.e.} $\alpha(j)=i$, $\alpha(i)=1$, 
$\alpha(1)=j$ and $\alpha(\ell)=\ell$ for any other index $\ell$. 
Then we have 
$$u_\alpha=\Bigl(u(j),...,\underbrace{u(1)}_{i},...,\underbrace{u(i)}_j,...,u(q)\Bigl).$$
Let also $\beta$ be the cycle $(1,j)$, {\em i.e.} $\beta$ exchanges $1$ with $j$ and keeps all other indexes unchanged. 
Then, as $u(i) \le u(j)$, $\beta$ is a reordering permutation of $u_{\alpha}$. 
\item if $i=1$, we let both $\alpha$ and $\beta$ be the cycle $(i,j)$. Then, 
$$u_\alpha=\Bigl(u(j),...,\underbrace{u(i)}_j,...,u(q)\Bigl),$$
so $\beta$ reorders $u_\alpha$. 
\item if $j=1$, we let $\alpha$ be the identity and $\beta$ be the cycle $(i,j)$. 
Then, 
$$u_\alpha=\Bigl(u(j),...,\underbrace{u(i)}_i,...,u(q)\Bigl),$$
and $\beta$ is again a reordering permutation of $u_\alpha$. 
\end{enumerate}
Now, fix $s$ and $t$ and let $\bar z$ be the following ordered element of $\R^q$,
$$\bar z=\left(t-s,t,...,t\right).$$
Let $F$ be a symmetric convex mapping~: $\R^q \to \R.$ 
In all cases, observe that $\alpha\circ\beta$ is the cycle $(1,i)$ (which amounts to the identity in case (2)), hence  
$$\left(u_{\alpha\circ\beta}-\bar z\right)(1)=u(i)+s-t.$$
Therefore, from the symmetry of $F$ and assertion (i) of Lemma \ref{lemma:schur-convexe1}, we have that
\begin{equation}
\label{eq:scx1}
F\left(u+s\gre_i-t\mathbf 1\right) = F\left(u_{\alpha\circ\beta} - \bar z\right).
\end{equation}
All the same, we have in all cases that 
$$
\left(u_{\alpha}-\bar z\right)(1) =u(j)+s-t$$
and all other coordinates of $u_{\alpha}-\bar z$ equal $u(\ell)-t$ for some index $\ell$. Therefore, by symmetry we have 
\begin{equation}
\label{eq:scx2}
F\left(u+s\gre_j-t\mathbf 1\right) = F\left(u_{\alpha} - \bar z\right).
\end{equation}
Now, as $\beta$ reorders $u_\alpha$, from claim (ii) of Lemma \ref{lemma:schur-convexe1} we have that 
$$u_{\alpha\circ\beta}-\bar z \lescx u_{\alpha}-\bar z,$$
which entails that $F\left(u_{\alpha\circ\beta}-\bar z\right) \le F\left(u_{\alpha}-\bar z\right).$ 
This together with (\ref{eq:scx1}) and (\ref{eq:scx2}) implies that 
$$F\left(u+s\gre_i-t\mathbf 1\right) \le F\left(u+s\gre_j-t\mathbf 1\right).$$ 
As this is true for any such mapping $F$, assertion (i) of Lemma \ref{lemma:schur-convexe1} concludes the proof. 
\end{proof}

Let us now introduce the following orderings,  

\begin{definition}
\label{def:lecx}
Let $q\ge 1$. For all $u,v \in \left(\R_+\right)^q$,
\begin{itemize}
\item we denote $u \lecx v$ whenever 
$$\sum_{i=k}^q \bar u(i) \le \sum_{i=k}^q \bar v(i),\mbox{ for all }k \in \llbracket 1,q \rrbracket.$$
\item For any $p \le q$, we write $u \prec_p^q v$ if 
\[\left\{\begin{array}{ll}
u &\lecx v\\
u(\ell) &\le v(\ell):\,\ell \in \llbracket p,q \rrbracket
\end{array}\right..\]
\end{itemize}
\end{definition}

\noindent The partial semi-ordering "$\lecx$" is nothing but a variant of "$\lescx$" to vectors of different total sums. 
First observe the following results, 
   
\begin{lemma}
\label{lemma:ordreétoile}
For all vectors $u,v \in \overline{\left(\R_+\right)^q}$ such that $u \lecx v$, 
\begin{itemize}
\item[(i)] for all $x\in\R$,  
$$\left[u-x\mathbf 1\right]^+ \lecx \left[v-x\mathbf 1\right]^+.$$
\item[(ii)] for all $j \in \llbracket 1,q \rrbracket$ such that $u(j) \le v(j)$, for all $y \in \R_+$, 
$$u+y\gre_j \lecx v+y\gre_j.$$
\end{itemize}
\end{lemma}
\begin{proof}
\begin{itemize}
\item[(i)] The result is trivial if $u(q) \le x$. Else, for all $k \in \llbracket 1,q \rrbracket$, for some $\ell \ge k$,  
\begin{align*}
\sum_{i=k}^S \left( u(i)-x\right)^+ = \sum_{i=\ell}^S (u(i)-x)
                                   \le \sum_{i=\ell}^S( v(i)-x)
                                   %&\le \sum_{i=\ell}^n (v(i)-x)^+\\
                                   \le \sum_{i=k}^S ( v(i)-x)^+.
\end{align*}
\item[(ii)] For all $k >j$, 
\begin{align*}
%\label{eq:shineon1}
\sum_{i=k}^S \left(\overline{u+y\gre_j}\right)(i) %&= u(k) \vee (u(j)+y) + \sum_{i=k+1}^n u(i)\\
                                &= \left(\sum_{i=k}^S  u(i) \right) \vee \left( u(j)+y+\sum_{i=k+1}^S  u(i)\right)\\
                                &\le \left(\sum_{i=k}^S  v(i) \right) \vee \left( v(j)+y+\sum_{i=k+1}^S  v(i)\right)\\
                                &=\sum_{i=k}^S \left(\overline{v+y\gre_j}\right)(i),                                     
\end{align*}
whereas for all $k \le j$, 
\begin{align*}
%\label{eq:shineon1}
\sum_{i=k}^S \left(\overline{u+y\gre_j}\right)(i) %&= u(k) \vee (u(j)+y) + \sum_{i=k+1}^n u(i)\\
                                     = \sum_{i=k}^S  u(i)+y 
                                     \le \sum_{i=k}^S  v(i) +y
                                     =\sum_{i=k}^S \left(\overline{v+y\gre_j}\right)(i).                                     
\end{align*}
\end{itemize}
\end{proof}

\noindent Recall the definition (\ref{eq:recurSCS0SP2}). We have the following result, establishing a comparison 
between the allocation policies JSW and 'Join the $p$-th Shortest Workload', for systems of the same size. 

\begin{proposition}
\label{pro:mainJSW}
Let $S\ge 1$ and $p \le S$.   
Then, for all $u,v \in \overline{\left(\R_+\right)^S}$ and all $s,t \ge 0$ we have that
$$u \prec_p^S v \,\Longrightarrow \, \Phi^{0,S}(u,s,t) \prec_p^S \Phi^{0,S,p}(v,s,t).$$
\end{proposition}
\begin{proof} Fix two non-negative numbers $s,t$ and two vectors 
$u,v \in \left(\R_+\right)^S$ such that $u \prec_p^S v$. 
First, we have for all $\ell \in \llbracket p,S \rrbracket$ that 
  \begin{align*}
    \Phi^{0,S}(u)(\ell) 
    &=\Bigl[u(\ell) \vee \bigl(\left(u(1)+s\right)
    \wedge u(\ell+1)\bigr) -t\Bigr]^+\\
%  & \le \left[u(\ell) \vee
%    \left(\left(u(P)+s\right)\wedge
%      u(\ell+1)\right) -t\right]^+\\
    & \le \left[v(\ell) \vee \left(\left(v(p)+s\right)\wedge 
    v(\ell+1)\right) -t\right]^+\\
    &=\Phi^{0,S,p}(v)(\ell).
  \end{align*}
This implies in particular that   
$$\sum_{i=k}^S \Phi^{0,S}(u)(i) \le \sum_{i=k}^S \Phi^{0,S,p}(v)(i),\,\mbox{ for all }k \ge p.$$ 
Consequently, %in order to verify (\ref{eq:condalloc1}) 
it suffices to check that for all $k < p$, 
\begin{equation}
  \label{eq:condallocoptimale3}
  \sum_{i=k}^S \Phi^{0,S}(u)(i)\le \sum_{i=k}^S \Phi^{0,S,p}(v)(i). 
\end{equation}
We have for all such $k$,
\[\left\{\begin{array}{ll}
  \displaystyle\sum_{i=k}^S \Phi^{0,S}(u)(i) &=\displaystyle\sum_{i=k+1}^{S} \left[u(i)-t\right]^+
  +\left[\left(u(1)+s\right)\vee u(k)-t\right]^+;\\
  \displaystyle\sum_{i=k}^S \Phi^{0,S,p}(v)(i) 
  &=\displaystyle\sum_{i=k;\,i\ne p}^{S} \left[v(i)-t\right]^+
  +\left[v(p)+s-t\right]^+.
\end{array}\right.\]
If $u(k) \ge u(1)+s$, then 
\begin{equation*}
  % \label{eq:condallocoptimale3}
  \sum_{i=k}^S \Phi^{0,S}(u)(i)=\sum_{i=k}^{S} \left[u(i)-t\right]^+ 
  \le \sum_{i=k}^{S} \left[v(i)-t\right]^+
  \le \sum_{i=k}^S \Phi^{0,S,p}(v)(i),
\end{equation*}
where we used item (i) of Lemma \ref{lemma:ordreétoile} in the first inequality. Consequently only the case where 
$u(k) < u(1)+s$ remains to be treated. 

For this, first observe that the vector $u$ is fully ordered, so we have $u(1) \le u(p)$. 
We can apply Lemma \ref{lemma:schur-convexe3} to $q=S-k+1$ and the vector 
$w=\left(u(1),u(k+1),...,u(S)\right)$: as the map 
\[\left\{\begin{array}{ll}
\left(\R_+\right)^{S-k+1} &\to \R\\
u &\mapsto \sum_{i=1}^{S-k+1} u^+
\end{array}\right.\]
is convex and symmetric, we have that 
\begin{align}
  \sum_{i=k}^S \Phi^{0,S}(u)(i)
  &=\left[u(1)+s-t\right]^++\sum_{i=k+1}^S \left[u(i)-t\right]^+\nonumber\\
  &\le \left[u(1)-t\right]^+
      +\sum_{i=k+1;i\ne p}^{S}\left[u(i)-t\right]^++\left[u(p)+s-t\right]^+\nonumber\\
  %\le \left[u(1)-t\right]^+
  %    +\sum_{i=k+1;i\ne S-N+1}^{S}\left[u(i)-t\right]^+\\
  %\shoveright{+\left[u(S-N+1)+s-t\right]^+}\\
  &\le \sum_{i=1;i\ne p}^{S}\left[u(i)-t\right]^+
     +\left[u(p)+s-t\right]^+.\label{eq:highhopes4}
  \end{align}
On the other hand, as $u(p) \le v(p)$, the assertions (ii) and then (i) of Lemma 
\ref{lemma:ordreétoile}, entail that
$$\left[u+s\gre_{p}-t\mathbf 1\right]^+ \lecx  
\left[v+s\gre_{p}-t\mathbf 1\right]^+.$$
In particular,
\begin{equation}
  \label{eq:highhopes5}
  \sum_{i=k;i\ne p}^{S}\left[u(i)-t\right]^++\left[u(p)+s-t\right]^+
  \le \sum_{i=k;i\ne p}^{S}\left[v(i)-t\right]^++\left[v(p)+s-t\right]^+,\, 
\end{equation}
and we deduce (\ref{eq:condallocoptimale3}) for $k<p$ such that $u(k) < u(1)+s$ 
from (\ref{eq:highhopes4}) and (\ref{eq:highhopes5}).
\end{proof}

\begin{proof}[Proof of Theorem \ref{thm:mainJSW}]
Consider two JSW systems such that $W_0^{0,S} \prec W_0^{0,S\prm}$, and fed by the same input $\suiten{\left(\sigma_n,\tau_n\right)}$.  
Let for all $n$, the following element of 
$\left(\R_+\right)^S$, 
\begin{equation}
\label{eq:defWtilde}
\tilde W_n^{0,S\pm} = \left(\underbrace{0,...,0}_{S-S\pm},W_n^{0,S\pm}\right).
\end{equation}
Observe that we have in $\left(\R_+\right)^S$, 
\begin{equation}
\label{fuckLMAC1}
W_0^{0,S} \prec_{S-S\prm+1}^S \tilde W_0^{0,S\pm}.
\end{equation}
Now, clearly, $\suite{\tilde W^{0,S\pm}}$ is nothing but the workload profile sequence of a 
'Join the $S-S\prm+1$-th shortest Workload' system of $S$ servers, of initial value $\tilde W^{0,S\pm}_0$. 
Therefore we have from (\ref{eq:recurJpSW}) that 
$$\tilde W_{n+1}^{0,S\prm}=\Phi^{0,S,S-S\prm+1}\left(\tilde W_{n}^{0,S\prm},\sigma_n,\tau_n\right),\,n\ge 0.$$
Consequently, (\ref{fuckLMAC1}) and Proposition \ref{pro:mainJSW} entail with a simple induction that for all $n\ge 0$, 
$W_n^{0,S} \prec_{S-S\prm+1}^S \tilde W_n^{0,S\pm}$, 
which is equivalent to saying that $W_n^{0,S} \prec W_n^{0,S\prm}$ for all $n$.   
\end{proof}

\section{At equilibrium}
\label{sec:steady}

Transient pathwise comparisons such as the results of Section \ref{sec:main} are particularly suitable for handling comparisons of systems at equilibrium 
using ergodic-theoretical representations. Let us describe in this way, how Theorem \ref{thm:mainJSW} can be transposed to the steady state of the system. 

\subsection*{Palm representation of the system}
Under the assumption that the input $\suiten{\left(\sigma_n,\tau_n\right)}$ is time-stationary and ergodic, we can adopt the so-called  
Palm representation of the queue, which we now quickly describe. The reader is referred to the classical monograph \cite{BacBre02} for a more complete picture. 

We work on the Palm space $\bar{\mathcal Q}=\left(\bar\Omega,\bar{\mathcal{F}}, \bar{\mathbf{P}}, \theta\right)$ of $\suiten{\left(\sigma_n,\tau_n\right)}$. 
We thereby assume that $\bar{\mathbf P}$-a.s. a customer enters the system at time 0 (the latter is denoted by $C_0$), requesting a service time of duration $\sigma$, and the following customer $C_1$ enters the system at time $\tau$. As a consequence of these basic assumptions, customer $C_n$ requests a service duration $\sigma\circ\theta^n$, and the inter-arrival epoch between the arrivals of $C_n$ and $C_{n+1}$ equals $\tau\circ\theta^n$. 
In particular, the sequences $\suiten{\left(\sigma_n,\tau_n\right)}$ and $\suiten{\left(\sigma\circ\theta^n,\tau\circ\theta^n\right)}$ have the same marginals, and the stationary ergodic assumptions imply that $\theta$ is $\mathbf {P}$-stationary and ergodic, that $\sigma$ and $\tau$ are integrable and $\tau>0$, $\mathbf P$-a.s..   

In these settings, a stochastic recursive sequence (SRS) is defined in the following way: if $X$ is a $E$-valued r.v. and $\Psi$ is a random mapping from $E$ 
to itself, the SRS starting from $X$ and driven by $\Psi$ is defined by
\[
\left\{\begin{array}{ll}
X_n&=X;\\
X_{n+1}&=\Psi\circ\theta^n\left(X_n\right),\,n\in\N.
\end{array}\right.\]
It is then routine to check that a time-stationary sequence having, on the original probability space $\mathcal Q$, 
the same distribution as $\{X_n\}$, corresponds to a solution $X$ defined on $\bar{\mathcal Q}$, to the functional equation  
\begin{equation}
\label{eq:recurstat}
X\circ\theta=\Psi\left(X\right),\bar{\mathbf P}-\mbox{ a.s..} 
\end{equation}
In the particular case where $E$ is a partially ordered (say, by "$\le$") Polish space having a minimal point $\mathbf 0$, if 
$\Psi$ is $\bar{\mathbf P}$-a.s. $\le$-nondecreasing and continuous, Loynes's Theorem (\cite{Loynes62}, see as well section 2.5.2 of \cite{BacBre02}) states 
that a solution $M_\infty$ to (\ref{eq:recurstat}) is given by the almost sure limit of the so-called {\em Loynes's sequence} $\suite{M}$, defined by  
\[\left\{\begin{array}{ll}
M_0&=0;\\
M_{n+1}&=\Psi\circ\theta^{-1}\left(M_n\circ\theta^{-1}\right),\,n\in\N. 
\end{array}\right.\] 
Moreover, the latter solution is {\em minimal}, in the sense that $M_\infty \le Y$ $\bar{\mathbf P}$-a.s. for any other 
solution $Y$ of (\ref{eq:recurstat}).

\subsection*{Join the Shortest Workload systems}
On the palm space $\bar{\mathcal Q}$ of arrivals of service, it is clear that the workload profile sequence of a JSW system of $S$ servers 
is a SRS driven by the map $\Psi^S:=\Phi^{0,S}(.,\sigma,\tau)$. 
Therefore a stationary service profile is a solution to the equation 
\begin{equation}
\label{eq:recurstatJSW}
V\circ\theta=\Psi^S(V),\,\bar{\mathbf P}-\mbox{ a.s..}
\end{equation}
Clearly, if "$\le$" denotes again the coordinate-wise ordering in $\R^S$, $\Psi^S$ is a.s. $\le$-non-decreasing and continuous, 
so Loynes's Theorem allows to construct the $\le$-minimal solution $W^{0,S}_\infty$ of (\ref{eq:recurstatJSW}). 
Moreover, all $S$ coordinates of $W^S$ are a.s. finite whenever $\esp{\sigma} < S\esp{\tau}$ (see \cite{Neu83}).  
We have the following result, 
\begin{theorem}
  \label{thm:mainsteady}
  Let $S\prm \le S$. 
  If $\esp{\sigma} < S\prm\esp{\tau}$, then the minimal solutions of (\ref{eq:recurstatJSW}) respectively for $S$ 
 and $S\prm$ servers, are such that 
$W^{0,S}_\infty \prec W^{0,S\prm}_\infty,\,\bar{\mathbf P}-\mbox{ a.s.,}$ 
in other words 
\begin{align}
W^{0,S}_\infty(S-i) &\le W^{0,S\prm}_\infty(S\pm-i),\,i\in \llbracket 0,S\prm-1\rrbracket;\label{eq:comparesteadyJSW1}\\
\mathbf W^{0,S}_\infty &\le \mathbf W^{0,S\prm}_\infty.\label{eq:comparesteadyJSW2}
\end{align}
\end{theorem}
\begin{proof}
Denote by $\{M_n^{S}\}$ (resp. $\{M_n^{S\prm}\}$), Loynes's sequence for the service profile of the system of $S$ 
(resp. $S\prm$) parallel queues.  
Let us consider for all $n\in\N$, the 
$\overline {\left(\R_+\right)^S}$-valued r.v. $\tilde M^{S\prm}_n$ defined likewise (\ref{eq:defWtilde}), {\em i.e.}  
$$\tilde M^{S\prm}_n=\left(\underbrace{0,\,...\,,0}_{S-S\prm\text{ terms }},M^{S\prm}_n(1),M^{S\prm}_n(2),...,M^{S\prm}_n(S\prm)\right).$$
It is easily checked that for all $n$, $\tilde M_{n+1}^{S\prm}=\Psi^{S-S\prm+1}\circ\theta^{-1}\left(\tilde M^{S\prm}_n\circ\theta^{-1}\right)$.  
% where $\Phi^{\sS-\sN+1}$ is defined by (\ref{eq:recurstatJSSW}) for $P:=S-N+1$.  
Therefore, as $\suite{M^{S}}$ and $\suite{\tilde M^{S\prm}}$ starts from the same initial value $\mathbf 0$, it follows again from a straightforward induction and Proposition \ref{pro:mainJSW}, that a.s. for all $n\in\N$, 
$ M^{S}_n \prec_{S-S\prm+1}^S \tilde M^{S\prm}_n$, or in other words $M^{S}_n \prec M^{S\prm}_n$. 
The result follows to the limit. 
\end{proof}
The latter result makes precise, in what sense the largest system performs better in steady state: from 
(\ref{eq:comparesteadyJSW2}) 
it minimizes the stationary total workload, and from (\ref{eq:comparesteadyJSW1}), the stationary offered waiting time since   
  \begin{equation}
  \label{eq:condallocoptimale3bis}
  W^{0,S}_\infty(1) \le W^{0,S\prm}_\infty(1),\,\bar{\mathbf P}-\mbox{ a.s..}
  \end{equation}

\begin{remark}
\label{remark:multiple}
Let us emphasize that these results can be adapted to multiple server queues operating in First Come, First Served: at equilibrium,   
Theorem \ref{thm:mainsteady} means that the server having the largest virtual workload among $S$ servers is almost surely less loaded than that of a system 
of $S\prm$ servers, and all the same for the second largest, the third largest, {\em etc...}, and that the total workload 
of the system of $S$ servers is a.s. less than that of the system of $S\prm$ servers. (\ref{eq:condallocoptimale3bis}) states that the system of $S$ servers 
offers almost surely a smaller waiting time than that of size $S\prm$. 
\end{remark}

\begin{remark}
\label{remark:Strassen}
All these pathwise comparison results can be translated from the Palm representation that is adopted here, onto the primary time-stationary queue, using a classical representation argument {\em \`a la} Strassen (see for instance Chapter 2 in \cite{BacBre02}). 
More precisely, it can be deduced from Theorem \ref{thm:mainsteady} that on the original 
probability space $\mathcal Q$, the time-stationary distributions of the service profiles $W^{S}$ and $W^{S\prm}$ satisfy 
\begin{align*}
W^{S}(S-\ell) &\le_{\text{st}}  W^{S\prm}(S\prm-\ell)\mbox{ for all }\ell \in \llbracket 0,S\prm-1 \rrbracket;\\
\mathbf W^{S} &\le_{\text{st}} \mathbf W^{S\prm},
\end{align*}
where "$\le_{\text{st}}$" denotes the strong (increasing) stochastic ordering associated to $\mathbf P$. 
\end{remark}

\section{Toward a pathwise comparison of allocation policies}
\label{sec:comparealloc}
To a certain extend, the latter results can be interpreted in terms of allocation policies for parallel queues. 
Fix the number $S$ of parallel queues, and let $p \le S$. Let $\suite{W^{0,S,p}}$ be the workload profile sequence of a system of $S$ servers applying the 'Join the $p$-th shortest workload' policy. As previously mentioned, the sequence $\suite{W^{0,S,p}}$ 
obeys the recursion 
$$W^{0,S,p}_{n+1}=\Phi^{0,S,p}\left(W^{0,S,p}_{n},\sigma_n,\tau_n\right),\,n\in\N,\,\mbox{ a.s..}$$ 
Denote again $\suite{W^{0,S}}$, the workload profile sequence of a system of $S$ servers under the JSW policy. 
Clearly, if the two systems are fed by the same input, from Proposition \ref{pro:mainJSW} we have that  
$$W^{0,S}_0 \prec_p^S W^{0,S,p}_0\,\Longrightarrow\,W^{0,S}_n \prec_p^S W^{0,S,p}_n,\mbox{ for all }n\ge 0.$$   
However, the latter allocation policy is irrelevant for applications: no more than $S-p+1$ servers will eventually admit customers, or in other words, at least $p-1$ servers will be left inactive from the first time in which they are idle. 
From that instant on, the system thus becomes a JSW system of $S-p+1$ servers. 

A maybe more interesting question is the following: consider now the allocation policy sending each incoming customer to the 
$p$-th smallest workload, {\em or to an empty queue}, if any. Let $\suite{\hat W^{0,S,p}}$ be the corresponding service profile sequence. 
It is clearly driven by the following map: for all $s,t$,  
\[\hat \Phi^{0,S,p}(.,s,t):\left\{\begin{array}{ll}
\overline {(\R_+)^S} &\to \overline {(\R_+)^S};\\
u &\mapsto \Phi^{0,S}\left(u,s,t\right)\ind_{\{u(1)=0\}}\,+\,\Phi^{0,S,p}\left(u,s,t\right)\ind_{\{u(1)>0\}}
\end{array}\right..\]
In order to compare the service profiles as above, we can aim at a pointwise comparison of the mappings 
$\Phi^{0,S,p}$ and $\hat\Phi^{0,S,p}$ for the ordering "$\prec_{p}^S$". However the situation is much more intricate in the present case: let $u$ and $v$, two elements of $\overline {(\R_+)^S}$ such that $u \prec_{p}^S v$. Then, 
it can be easily checked that for any $s,t$, $\Phi^{0,S}(u,s,t) \prec_{p}^S \hat\Phi^{0,S,p}(v,s,t)$ if $v(1)>0$ 
(this is Proposition \ref{pro:mainJSW}), or $u(1)=0$ (see item (ii) of Lemma 5 of \cite{Moy16}). 
But one can construct examples for which $\Phi^{0,S}(u,s,t) \not\prec_{p}^S \hat\Phi^{0,S,p}(v,s,t)$ with positive probability if $u(1)>0$ and $v(1)=0$.  
 
Therefore there is no pathwise ordering of the transient profile sequences in this context 
(at least for the ordering "$\prec_{p}^S$"). We conjecture, however, that the latter holds true at equilibrium. Advances on this open problem can be found in \cite{Moy16}, together with the stability 
study of this alternative system.

\providecommand{\bysame}{\leavevmode\hbox to3em{\hrulefill}\thinspace}
\providecommand{\MR}{\relax\ifhmode\unskipPace\fi MR }
\providecommand{\MRhref}[2]{%
  \href{http://www.ams.org/mathscinet-getitem?mr=#1}{#2}
}
\providecommand{\href}[2]{#2}

\end{document}